\documentclass[conference]{IEEEtran}
\usepackage{amsmath,amssymb,amsthm,epsfig,color,graphicx}
\usepackage{enumerate,url,algorithm,algorithmic,verbatim,balance}

\newcommand \bzero{\mathbf{0}}
\newcommand \bone{\mathbf{1}}

\newcommand \bb{\mathbf{b}}

\newcommand \bg{\mathbf{g}}

\newcommand \bp{\mathbf{p}}
\newcommand \bq{\mathbf{q}}

\newcommand \bu{\mathbf{u}}
\newcommand \bv{\mathbf{v}}
\newcommand \bw{\mathbf{w}}
\newcommand \bx{\mathbf{x}}
\newcommand \by{\mathbf{y}}
\newcommand \bz{\mathbf{z}}

\newcommand \bR{\mathbf{R}}

\newcommand \bX{\mathbf{X}}

\newcommand \btheta{\boldsymbol{\theta}}

\newcommand \blambda{\boldsymbol{\lambda}}

\newcommand \bpi{\boldsymbol{\pi}}

\newcommand \mcQ{\mathcal{Q}}

\begin{document}
	
	\title{Deep Learning for Reactive Power Control of Smart Inverters under Communication Constraints}
	
	\author{
	\IEEEauthorblockN{Sarthak Gupta, Vassilis Kekatos, and Ming Jin}
	
	\IEEEauthorblockA{
	\textit{Bradley Dept. of Electrical \& Computer Engineering}\\
	\textit{Virginia Tech}, Blacksburg, VA 24061, USA\\
	Emails: \{gsarthak,kekatos,jinming\}@vt.edu}
	}
	
\maketitle
	
\begin{abstract}
Aiming for the median solution between cyber-intensive optimal power flow (OPF) solutions and subpar local control, this work advocates deciding inverter injection setpoints using deep neural networks (DNNs). Instead of fitting OPF solutions in a black-box manner, inverter DNNs are naturally integrated with the feeder model and trained to minimize a grid-wide objective subject to inverter and network constraints enforced on the average over uncertain grid conditions. Learning occurs in a quasi-stationary fashion and is posed as a stochastic OPF, handled via stochastic primal-dual updates acting on grid data scenarios. Although trained as a whole, the proposed DNN is operated in a master-slave architecture. Its master part is run at the utility to output a condensed control signal broadcast to all inverters. Its slave parts are implemented by inverters and are driven by the utility signal along with local inverter readings. This novel DNN structure uniquely addresses the small-big data conundrum where utilities collect detailed smart meter readings yet on an hourly basis, while in real time inverters should be driven by local inputs and minimal utility coordination to save on communication. Numerical tests corroborate the efficacy of this physics-aware DNN-based inverter solution over an optimal control policy.
\end{abstract}
	

\section{Introduction}\label{sec:intro}
\renewcommand{\thefootnote}{\fnsymbol{footnote}}
\footnotetext{This work was supported in part by the NSF-1751085 grant.} 
\allowdisplaybreaks

Distribution grids are currently challenged by voltage fluctuations due to the proliferation of distributed energy resources (DERs). The voltages experienced at buses of a feeder depend heavily on the power injected or withdrawn, while the power generated by a PV under intermittent cloud coverage may vary by 80\% within one-minute intervals~\cite{Turitsyn11}. The inverters interfacing DERs have been suggested as a promising fast-responding mechanism and are now allowed to provide reactive power support per the amended IEEE~1547 standard. If properly orchestrated, inverters can regulate nodal voltages and/or reduce ohmic line losses. Nonetheless, coordinating hundreds of inverters in real-time is a formidable task. 

The literature on inverter control can be broadly classified into optimization- and learning-based approaches. The former class includes approaches where inverter control is posed as an optimal power flow (OPF) problem. Under a centralized OPF setup~\cite{FCL},~\cite{ergodic}, the utility reads the values of solar generation and loads, solves an OPF, and communicates the optimal setpoints to inverters. To avoid any cyber overhead, inverter setpoints can be decided using simple Volt/VAR or Watt/VAR control rules driven by local readings~\cite{Turitsyn11}. Nonetheless, the equilibria of such rules do not coincide with the sought OPF solutions and can be subpar~\cite{KZGB15},~\cite{Jabr18}.

Learning-based approaches shift the computational effort offline, and perform numerically less intensive tasks during real-time operation. Learning-based approaches can be further clustered into the \emph{OPF-then-learn} and the \emph{OPF-and-learn} philosophies. According to the former, one first solves a large number of OPF instances parameterized by their inputs (solar/load conditions). The pairs of OPF inputs or instances and OPF minimizers are subsequently used for the ML model to learn the OPF mapping in a supervised manner. In real time, the ML model approximates OPF decisions on the fly as soon as it is presented with a new OPF instance. Under this paradigm, references \cite{Dobbe19} and \cite{Kara18} use kernel--based regression to learn inverter control rules. DNNs have alternatively been employed to learn OPF solutions under a linearized~\cite{DeepOPFPan19}; or an exact AC grid model~\cite{ZamzamBaker19}, \cite{GuhaACOPF}, \cite{ZhangDNNOPF}, \cite{RibeiroGNNOPF}.

Rather than fitting OPF minimizers, the \emph{OPF-and-learn} paradigm trains an ML model directly through an OPF in a single step. Therefore, it does not require solving multiple OPFs to generate a labeled training set. Under the \emph{OPF-and-learn} paradigm, reference \cite{JKGD19} adopts kernel-based learning to design inverter control rules, adjusted to grid conditions in a quasi-stationary fashion. Although rules can be learned using a convex program, the kernel functions have to be specified beforehand. In~\cite{Yang19}, inverter control rules are optimized along with capacitor status on/off decisions to minimize voltage deviations using a two-timescale reinforcement learning (RL) approach. Nonetheless, no feeder-level constraints are involved. Enforcing network constraints is challenging for learning-based OPF methods. One could heuristically project the ML prediction for the OPF solution~\cite{ZamzamBaker19}, \cite{JKGD19}. Other approaches to coping with constraints include penalizing constraint deviations~\cite{Kara18}, \cite{DeepOPFPan19}, \cite{JKGD19}, \cite{Yang19}; or enforcing constraints in a discounted sense~\cite{ZhangISU}. Reference~\cite{ZhangISU} models inverter policies as DNNs. It successively linearizes feeder constraints and updates policies continuously through communication exchanges between interconnected microgrids. A similar \emph{safe RL} learning scheme is put forth in~\cite{Nanpeng19}, but with a centralized implementation. 

A key promise of designing policies is to alleviate the cyber burden of inverter control. This critical aspect has been largely overlooked by the existing literature. In particular, references \cite{ZhangISU} and \cite{Nanpeng19}, which are most closely related to this work, update policies continuously and require considerable amounts of data to be communicated in real time. To account for this aspect, the contributions of this work are in two fronts: First, inverter policies are modeled as DNNs that are jointly trained in a quasi-stationary fashion, while feeder constraints are enforced explicitly in a stochastic sense. Second, a carefully designed DNN architecture accommodates application scenarios where inverter rules are driven by local measurements as well as a low-bandwidth control signal broadcast by the utility.

The rest of this work is organized as follows. Section~\ref{sec:problem} formulates the task of designing inverter control policies after reviewing an approximate grid model. Section~\ref{sec:learning} adopts a stochastic primal-dual algorithm to find the optimal inverter control policies. Section~\ref{sec:architecture} puts forth the novel communication-cognizant DNN-based inverter control architecture. The proposed schemes are evaluated using real-world solar generation and load data on the IEEE 13-bus feeder in Section~\ref{sec:tests}. Conclusions along with ongoing and future research directions are discussed in Section~\ref{sec:conclusions}.

\emph{Notation:} lower- (upper-) case boldface letters denote column vectors (matrices), and calligraphic symbols are reserved for sets. Symbol $^{\top}$ stands for transposition and $\|\bx\|_2$ denotes the $\ell_2$-norm of $\bx$. Vectors $\bzero$ and $\bone$ are respectively the vectors of all zeros and ones of appropriate dimensions.

\section{Grid Modeling and Problem Formulation}\label{sec:problem}
Consider a feeder with $N+1$ buses, including the substation indexed by $0$. Let $p_n+jq_n$ be the complex power injection at bus $n$. Its active power component can be decomposed as $p_n=p_n^g-p_n^c$, where $p_n^g$ is the solar generation and $p_n^c$ the inelastic load at bus $n$. Its reactive power component can be similarly expressed as $q_n=q_n^g-q_n^c$. If vectors $(\bp,\bq)$ collect the power injections at all non-substation buses, they can be decomposed as
$\bp=\bp^g-\bp^c~~\textrm{and}~~\bq=\bq^g-\bq^c$.
We refer to the values of (re)active loads and active solar generation at all non-substation buses as \emph{grid conditions}
\begin{equation}\label{eq:cond}
\bz:=[(\bp^c)^\top~~(\bq^c)^\top~~(\bp^g)^\top]^\top.
\end{equation}
Given $\bz$, the task of reactive power control by DERs aims at optimally setting $\bq^g$ to minimize a feeder-wide objective while complying with network and inverter limitations. Starting with the latter, the reactive power injected by inverter $n$ is limited by a given $\bar{q}_n^g$ due to apparent power limits. Apparent power constraints are local and will be collectively denoted by 
\begin{equation}\label{eq:Q}
\bq^g\in\mathcal{Q}:=\left\{\bq:|q_n^g|\leq \bar{q}_n^g\quad\forall n\right\}.
\end{equation}

Regarding feeder constraints, the focus is on confining voltages within the regulation range of $[0.97,1.03]$~per unit (pu). Albeit voltages are nonlinearly related to power injections, for simplicity we adopt a widely used linearized grid model~\cite{TJVT20}. According to this model, the vector of voltage magnitudes at all $N$ buses is approximately
\begin{equation}\label{eq:ldf}
\bv=\bR\bp+\bX\bq+v_0\bone    
\end{equation}
where $v_0$ is the substation voltage, while the symmetric positive semidefinite matrices $(\bR,\bX)$ depend on the feeder and are assumed to be known. If each voltage $v_n$ is to be maintained within $[\underline{v}_n,\overline{v}_n]$, the reactive power injections $\bq^g$ should satisfy the network constraints
\begin{equation}\label{eq:volt}
\bg(\bq^g,\bz):=\left[\begin{array}{c}
\bX\bq^g +\by-\overline{\bv}\\
-\bX\bq^g -\by+\underline{\bv}
\end{array}\right]\leq \bzero
\end{equation}
where vector $\by:=\bR(\bp^g-\bp^c)-\bX\bq^c+v_0\bone$ depends on $\bz$, and vectors $(\underline{\bv},\overline{\bv})$ contain the limits $(\underline{v}_n,\overline{v}_n)$ across buses.

According to the same grid model, ohmic losses on lines can be approximated as a convex quadratic function of power injections as $\bp^\top\bR\bp+\bq^\top\bR\bq$; see~\cite{TJVT20} for details. Upon defining $\bb:=2\bR\bq^c$, the part of ohmic losses that is dependent on the control variable $\bq^g$ can be approximated as
\begin{equation}\label{eq:losses}
\ell(\bq^g,\bz)=(\bq^g)^\top\bR\bq^g-\bb^\top\bq^g.
\end{equation}

We henceforth abuse notation and use $\bq$ in lieu of $\bq^g$. This should not cause any confusion since $\bq^c$ has been included in $\bz$. DER reactive setpoints $\bq$ can be found as the minimizer of
\begin{align}\label{eq:opf}
\min_{\bq\in\mcQ}~&~\ell(\bq,\bz)\\
\mathrm{s.to}~&~\bg(\bq,\bz)\leq \bzero.\nonumber
\end{align}
Under the linearized grid model, the approximate OPF task of \eqref{eq:opf} is a convex quadratic program (QP). Solving \eqref{eq:opf} can be computationally and communication-wise taxing if $\bz$ changes frequently. Moreover, by the time \eqref{eq:opf} is solved and decisions are downloaded to DERs, grid conditions $\bz$ may have changed rendering the computed setpoints obsolete. 

To account for the uncertainty in $\bz$, one may pursue a stochastic formulation such as~\cite{ergodic}
\begin{align}\label{eq:opf2}
\min_{\bq\in\mcQ}~&~\mathbb{E}[\ell(\bq,\bz)]\\ \mathrm{s.to}~&~\mathbb{E}[\bg(\bq,\bz)]\leq \bzero\nonumber
\end{align}
where the expectation $\mathbb{E}$ is with respect to $\bz$. Nonetheless, the obtained `one-size-fits-all' $\bq$ does not adapt to different $\bz$'s.

To come up with DER setpoints that are responsive to grid conditions, we resort to \emph{control policies or rules}, where the reactive power setpoint for each inverter $n$ is captured by a function $\pi_n(\bw_n;\btheta_n)$ acting upon a control input $\bw_n$ and is parameterized by vector $\btheta_n$. Ideally, inverter control policies should be driven by the complete $\bz$, that is $\bw_n=\bz$ for all $n$. Nevertheless, that would entail high communication overhead. If the utility knows the complete $\bz$, it might as well solve \eqref{eq:opf} and communicate the optimal setpoints to inverters. For an inverter control scheme to be communication-cognizant, the inputs $\bw_n$ should primarily involve local readings of $\bz$, such as $(p_n^g,q_n^g,p_n^g)$, and possibly few remote entries. Regarding the parameter vectors $\btheta_n$'s, these may be unique per inverter or share some entries as detailed in Section~\ref{sec:architecture}. To capture the aforementioned scenarios, let us abstractly refer to the vector of inverter policies $\pi_n(\bw_n;\btheta_n)$'s as 
\begin{equation}\label{eq:rule}
\bq(\bw)=\bpi(\bw;\btheta)
\end{equation}
where $\bw$ is the union of $\bw_n$'s and $\btheta$ the union of $\btheta_n$'s. 

The control policies for DERs can be found jointly by solving the constrained stochastic minimization
\begin{align}\label{eq:opf3}
P^*:=\min_{\btheta:\bpi(\bw;\btheta)\in\mcQ}~&~\mathbb{E}[\ell(\bpi(\bw;\btheta),\bz)]\\
\mathrm{s.to}~&~\mathbb{E}[\bg(\bpi(\bw;\btheta),\bz)]\leq \bzero\nonumber
\end{align}
over the parameter vector $\btheta$. Problem~\eqref{eq:opf3} couples policies in two ways. First, for a fixed $\bz$, policies are coupled across inverters through the cost and constraint functions since the entries of $\bq^g$ appearing in \eqref{eq:ldf} and \eqref{eq:losses} are now computed via \eqref{eq:rule}. Second, the expectations in \eqref{eq:opf2} and \eqref{eq:opf3} couple system's performance across OPF instances characterized by $\bz$. 

Local and linear policies of the form $\pi_n(\bw_n;\btheta_n)=\btheta_n^\top\bw_n$ have been previously studied for inverter control~\cite{Jabr18}, \cite{LinBitar18}, \cite{Baker18}. Nonetheless, the optimal policies $q_n(\bw_n)$ are not necessarily affine in $\bw_n$, especially when $\bw_n$ is only a partial observation of $\bz$. The grand challenge towards scalable inverter control is to design \emph{nonlinear control curves}. In~\cite{JKGD19}, we dealt with by modeling each $q_n(\bw_n)$ as a kernel-based support vector machine (SVM), and designing all rules jointly under an OPF formulation. The advantage of SVM-based policies is that they can be trained to optimality using convex optimization. Nonetheless, selecting the appropriate kernel and control inputs $\bw_n$'s can be challenging. Inspired by their field-changing performance in various engineering tasks, here we propose modeling inverter rules using DNNs, and train the parameters $\btheta$ in a data-driven physics-aware fashion.

\section{Primal-dual DNN Learning}\label{sec:learning}
Solving \eqref{eq:opf3} is challenging since it is a constrained stochastic minimization over a DNN. To train the inverter policy DNN, we adopt the stochastic primal-dual updates of~\cite{Ribeiro19}, which are briefly reviewed next. Consider the Lagrangian function of \eqref{eq:opf3}
\begin{align}\label{eq:lagrangian}
L(\btheta;\blambda)=  \mathbb{E}[\ell(\bpi(\bw;\btheta),\bz)] + \blambda^\top\mathbb{E}[\bg(\bpi(\bw;\btheta),\bz)]
\end{align}
where $\blambda$ is the vector of Lagrange multipliers corresponding to constraint \eqref{eq:opf3}. The dual problem can be posed as
\begin{align}\label{eq:dual}
D^*=\max_{\blambda\geq \bzero}\min_{\btheta:\bpi(\bw;\btheta)\in\mcQ} L(\btheta;\blambda).
\end{align}
Standard duality results predicate that $D^*\leq P^*$. When the primal problem is convex, the previous inequality typically holds with equality. Problem \eqref{eq:opf3} however is non-convex even if \eqref{eq:opf} is a convex QP, since the DNN mapping $\bpi(\bw;\btheta)$ is generally non-convex in $\btheta$. Nonetheless~\cite{Ribeiro19} establishes that: \emph{i)} under relatively mild conditions satisfied by \eqref{eq:opf3}, and \emph{ii)} if the underlying DNN architecture is rich enough, the duality gap $P^*-D^*$ is sufficiently small. This motivates solving \eqref{eq:opf3} through the primal-dual updates indexed by $k$~\cite{Ribeiro19}
\begin{subequations}\label{eq:pdupdate}
\begin{align}
    \btheta^{k+1}&=\big[\btheta^{k}-\mu_{\theta}\nabla_{\btheta}L(\btheta^{k};\blambda^{k})\big]_{\mcQ}\label{eq:pdupdate:p}\\
    \blambda^{k+1}&=\big[\blambda^{k}+\mu_{\lambda}\nabla_{\blambda}L(\btheta^{k+1};\blambda^{k})\big]_+\label{eq:pdupdate:d}
\end{align}
\end{subequations}
where the operator $[\cdot]_{\mcQ}$ projects $\btheta^{k+1}$ such that $\bpi(\bw;\btheta^{k+1})\in\mcQ$ for all $\bw$; operator $[\cdot]_+$ ensures $\blambda\geq\mathbf{0}$ at all times; and $(\mu_{\theta},\mu_{\lambda})$ are positive step sizes. Regarding $[\cdot]_{\mcQ}$, the DNN output corresponding to $q_n^g$ can be constrained within $[-\bar{q}_n^g,+\bar{q}_n^g]$ by using $\tanh(\cdot)$ as the output activation function and then scaling by the constant $\overline{q}_n^g$.

The updates in \eqref{eq:pdupdate} are complicated by the expectation operator. The probability distribution function (pdf) of $\bz$ (and hence $\bw$) may not be known beforehand. Even if it is known, propagating that pdf through nonlinear functions such as $\bpi(\bw;\btheta),\bz)$ is non-trivial. To deal with this, the primal-dual updates of \eqref{eq:pdupdate} can be surrogated by their stochastic approximation counterparts relying on samples of grid conditions. In particular, the utility is assumed to have a set of scenarios $(\bz^k,\bw^k)$ indexed by $k=1,\ldots,K$, with which the ensemble averages of \eqref{eq:lagrangian} are approximated as $\mathbb{E}[\ell(\bpi(\bw;\btheta),\bz)]\simeq \frac{1}{K}\sum_{k=1}^K\ell(\bpi(\bw^k;\btheta),\bz^k)$. To simplify the updates of \eqref{eq:pdupdate}, the sample averages can be approximated by a single scenario per iteration to yield the \emph{stochastic} primal-dual updates~\cite{Ribeiro19}
\begin{subequations}\label{eq:spdupdate}
\begin{align}
    \btheta^{k+1}&=\bigg[\btheta^{k}-\mu_{\theta}\big(\nabla_{\btheta}\ell^k- (\blambda^k)^\top\nabla_{\btheta}\bg^k\big)\bigg]_{\mcQ}\label{eq:spdupdate:p}\\
    \blambda^{k+1}&=\left[\blambda^{k}+\mu_{\lambda}\bg\left(\bpi(\bw^k;\btheta^{k+1}),\bz^k\right)\right]_+.\label{eq:spdupdate:d}
\end{align}
\end{subequations}
Here $\nabla_{\btheta}\ell^k$ is the gradient of $\ell(\bpi(\bw;\btheta),\bz)$ and $\nabla_{\btheta}\bg^k$ the Jacobian matrix of $\bg(\bpi(\bw;\btheta),\bz)$, both with respect to $\btheta$ and evaluated at $(\bw^k,\btheta^k,\bz^k)$. The updates are known to converge to a stationary point of \eqref{eq:opf3} for sufficiently small step sizes. 

For the objective and constraint functions of \eqref{eq:volt}--\eqref{eq:losses}, the needed sensitivities can be computed as
\begin{align*}
&\nabla_{\btheta}\ell^k=\left(\nabla_{\btheta}\bpi(\bw^k;\btheta^k)\right)^\top\left(2\bR\bpi(\bw^k;\btheta^k) - {\bb^k}\right)\\
&\nabla_{\btheta}\bg^k=\left[\bX~~-\bX\right]^\top\nabla_{\btheta}\bpi(\bw^k;\btheta^k).
\end{align*}
Here $\bb^k:=2\bR(\bq^c)^k$ and $\nabla_{\btheta}\bpi(\bw^k;\btheta^k)$ is the Jacobian matrix of the DNN output with respect to its weight parameters. The latter can be evaluated using gradient back-propagation across the DNN, a standard tool readily available in all DNN-related software. If the number of available grid scenarios $K$ is relatively small, additional scenarios can be synthesized by applying small perturbations on the available $\bz^k$'s. As customary in DNN training, the updates \eqref{eq:spdupdate} can be iterated over multiple epochs or in mini-batch forms.

It is worth contrasting the DNN input $\bw$ and the vector of grid conditions $\bz$. Despite some possible overlap, the two vectors are used differently. The former one feeds the DNN to compute the setpoints $\bq(\bw)=\bpi(\bw;\btheta)$. The latter one is involved in the OPF objective and constraint functions, i.e., it appears in $\bb$ for computing $\nabla_{\btheta}\ell$ and when evaluating $\bg(\bpi(\bw;\btheta),\bz)$. While $\bz$ should be known to the utility during training to perform the updates of \eqref{eq:spdupdate}, it is not needed during real-time operation. This resonates with the small/big data setup, since a utility has offline access to an extensive smart meter dataset of $\bz$'s; yet its control center and each inverter individually are driven by limited real-time data feeds. The updates of \eqref{eq:spdupdate} apply for inverters DNNs of arbitrary architecture. We next particularize the structure of $\bpi(\bw;\btheta)$ to comply with communication limitations in inverter control. 



\section{Communication-Cognizant DNN Architecture}\label{sec:architecture}

\begin{figure}[t]
    \centering
    \includegraphics[scale=0.28]{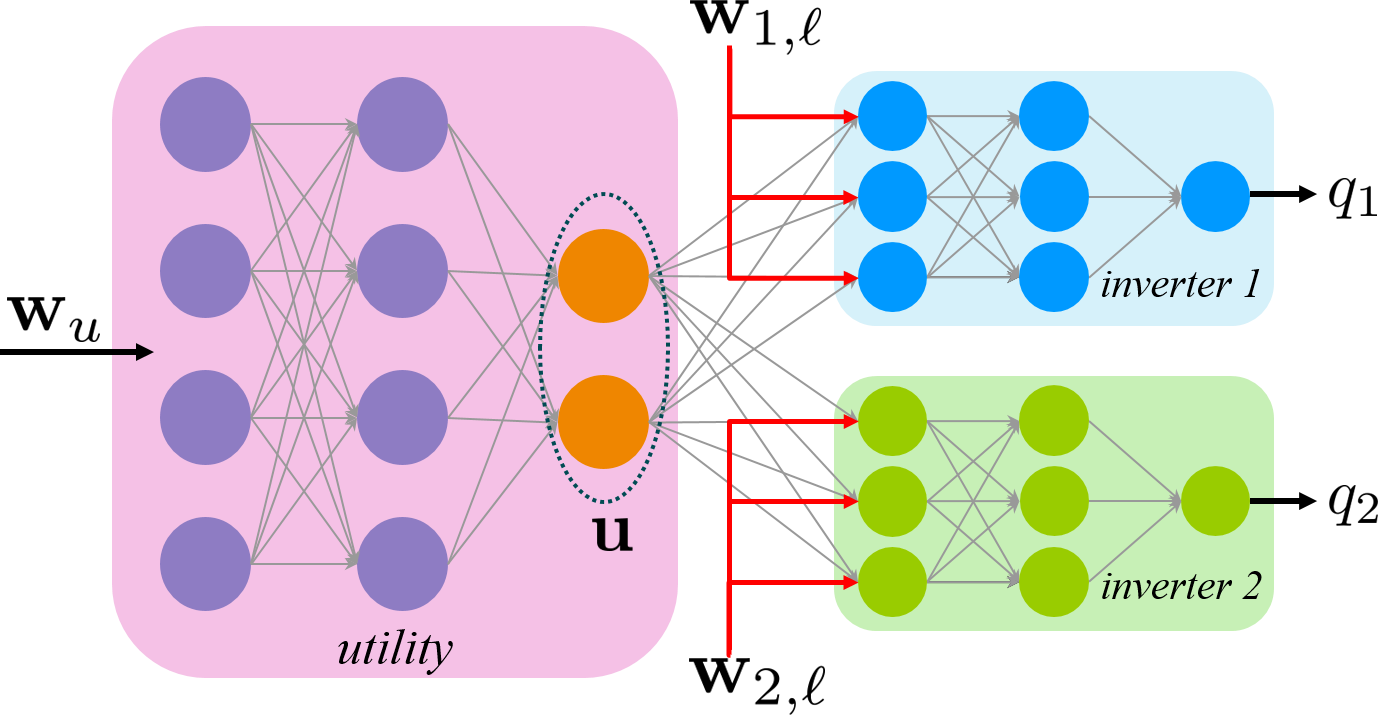}\\
    \includegraphics[scale=0.48]{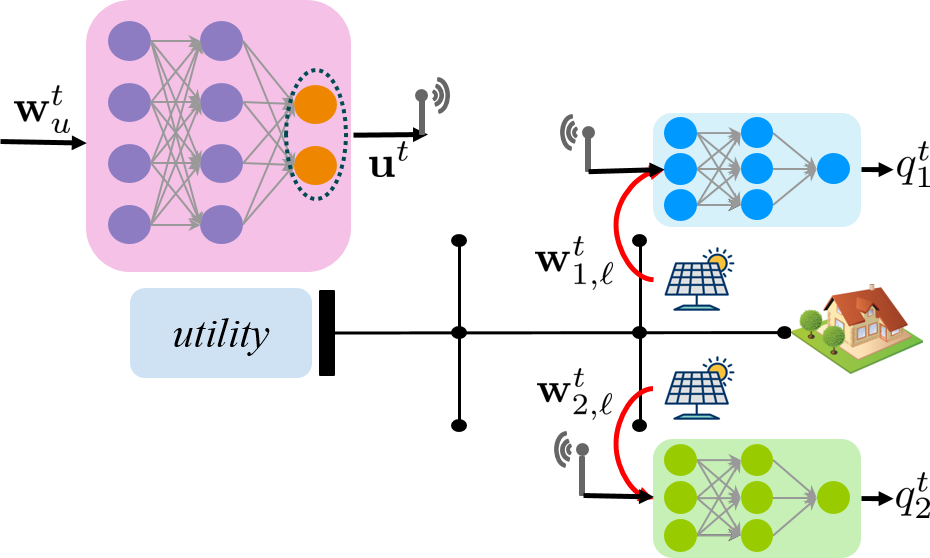}
    \caption{\emph{Top}: DNN $\bpi(\bw;\btheta)$ is organized in one utility sub-NN and inverter sub-NNs, all trained as a single DNN by the utility offline. \emph{Bottom}: During real-time operation, the utility sub-NN uses real-time data to compute and broadcast the control signal, while inverter sub-NNs are run at inverters.}
    \label{fig:training}
\end{figure}

To coordinate inverters on a tight communication budget, our proposed inverter policy DNN $\bpi(\bw;\btheta)$ comes with the two-tier architecture shown on Figure~\ref{fig:training} (top). Its first layers constitute the \emph{utility sub-NN}, while the final layers constitute the \emph{inverter sub-NNs}, one for each inverter. Figure~\ref{fig:training} shows only two inverters for simplicity. The utility sub-NN (shown in purple) is fully connected, is driven by input $\bw_u$, and outputs control $\bu$. Inverter sub-NNs (in blue and green) are disconnected from each other and both fed with the common control $\bu$. Each inverter sub-NN is also fed with its own local data $\bw_{n,\ell}$. The $n$-th inverter sub-NN predicts the setpoint $q_n$.

Inverter policy $n$ can be expressed as $q_n(\bw_n)=\pi_n(\bw_n;\btheta_n)$ where $\bw_n=[\bw_u^\top~\bw_{n,\ell}^\top]^\top$ and $\btheta_n$ collects the DNN parameters for the shared utility sub-NN and inverter sub-NN $n$. Vectors $\bw_{n,\ell}$ may carry local load and solar generation available on bus $n$. Input $\bw_u$ carries information available to the utility control center in real time. Such information can be power flow readings from major distribution lines, transformers, and/or voltage regulators. Vector $\bw_u$ may also carry the solar generation from a solar farm or any other DER that is telemetered in real time. Rather than actual grid measurements, vector $\bw_u$ may also include predictions the utility can make on grid conditions. For example, that could be the case if the utility uses cameras to monitor cloud coverage as a proxy to solar generation or temperature/humidity readings to load. 

During training and given grid scenario $\bz^k$, the inputs $\bw_u^k$ and $\bw_n^k$'s can be: \emph{i)} found readily as partial entries of $\bz^k$ (loads and solar generation); \emph{ii)} inferred from $\bz^k$ (a line flow can be computed through the power flow equations, or approximated as the sum of all downstream power injections); or \emph{iii)} found through historical data (dataset combining cloud coverage with solar generation). The particular structure of the proposed DNN with individualized inputs and partially connected layers can be easily implemented by skipping and masking connections, respectively.

\begin{algorithm}[t]
	\caption{Inverter control through DNN-based policies} \label{alg:1}
	\textbf{\emph{Training (utility side)}}
	\begin{algorithmic}[1]
		\STATE Collect grid scenarios $\{\bz^k\}_{k=1}^K$ from smart meter data
		\STATE Collect or calculate DNN inputs $\{\bw^k\}_{k=1}^K$
		\STATE Initialize $\btheta^0$ and $\blambda^0$
		\FOR{all $K$ scenarios and $E$ epochs}
		\STATE Update $\btheta$ using \eqref{eq:spdupdate:p}
		\STATE Update $\blambda$ using \eqref{eq:spdupdate:d}
		\ENDFOR
	    \STATE Download $\btheta$ parameters to inverter sub-NNs
	\end{algorithmic}
	\textbf{\emph{Real-time operation (utility \& inverter sides)}}
	\begin{algorithmic}[1]
	    \FOR{$t=0,1,\ldots,T,$}
		\STATE Utility receives $\bw_u^t$ from real-time telemetry
		\STATE Feed $\bw_u^t$ to utility sub-NN to compute $\bu^t$
		\STATE Utility broadcasts $\bu^t$ to inverters
		\FOR{each inverter $n$}
		  \STATE Inverter $n$ reads $\bu^t$ and local data $\bw_n^t$
    	  \STATE Feed $(\bu^t,\bw_n^t)$ to inverter sub-NN to decide $q_n^t$
		\ENDFOR
		\ENDFOR
	\end{algorithmic}
\end{algorithm}

Although trained as a whole, the inverter policy DNN $\bpi(\bw;\btheta)$ is implemented in parts; see bottom panel of Fig.~\ref{fig:training}. After training is completed for the upcoming 30- or 60~min period, the weights corresponding to inverter NNs are downloaded to inverters. A unique component of our DNN architecture is the control signal $\bu$, which is broadcast from the utility NN to inverter NNs. To save on downlink (utility to inverters) communications, signal $\bu$ is designed to be much shorter than $\bw_u$. Considering that $\bu$ is actually designed along with the operation of inverter sub-NNs through the OPF of \eqref{eq:opf3}, this signal carries all the information the utility can provide to coordinate inverters in a condensed form. Its broadcast nature further contributes to communication savings. The steps involved during the training and real-time operation of the proposed DNN are summarized in Algorithm~\ref{alg:1}. 

This DNN architecture can cater to a wide range of communication specifications. If no downlink communication is allowed in real time, the utility sub-NN can be ignored all together and inverter sub-NNs are driven based on local inputs. If downlink bandwidth is abundant, inverter sub-NNs can be dropped and inverter setpoints can be decided by the utility sub-NN in real time. Practical application scenarios are expected to lie somewhere between these two extremes, whence the hybrid architecture of Fig.~\ref{fig:training} becomes relevant.

\section{Numerical Tests}\label{sec:tests}
The proposed DNN-based inverter control was evaluated on a single-phase version of the IEEE 13-bus feeder. Real-world active load data was extracted for March 1, 2018, on a one-minute resolution from Pecan Street. Solar generation data was also added to buses $\{1,5,9,10,11,12\}$, out of which buses $\{9,12\}$ were equipped with inverters. Load time series were scaled so that monthly peaks were $7.5$ times the benchmark values. The same ratio was used to scale solar. Reactive loads were added with lagging power factors sampled from a uniform distribution between 0.9 and 1. The utility was assumed to have telemetry $\bw_u$ for the active line flows feeding buses $\{2,3,7\}$ from their parent buses.

\begin{figure}
    \centering
    \includegraphics[scale=0.4]{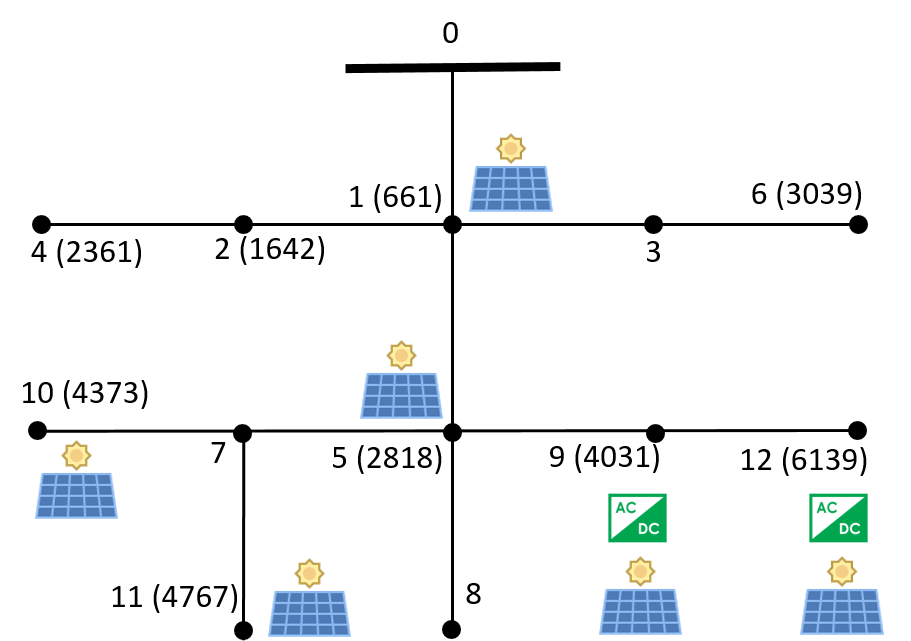}
    \caption{The IEEE 13-bus feeder. Numbers in parentheses indicate the house index from the Pecan Street dataset mapped to each bus.}
    \label{fig:IEEE13}
\end{figure}


The utility sub-NN was constructed using an input layer of dimension 3 and an output layer $\bu$ of dimension 1.Inverter sub-NN were made up of one input, hidden, and output layers of dimensions 5, 6 and 1, respectively. The local readings $\{p_n,q^c_n\}$ along with $\bu$ were fed as inputs to each inverter sub-NN $n$. Initial values for DNN parameters were uniformly sampled from the range $[-0.1, 0.1]$ and were updated using Adam with a learning rate of 0.01. The dual variables were all initialized at 0 and were updated with step sizes of 1 that decayed with the square-root of the iteration index~\cite{LKMG17}. Our approach was contrasted with an optimal policy $\bq(\bw)$ that directly solves \eqref{eq:opf2} without being confined to any DNN or other parameterization using dual decomposition~\cite{LKMG17}.


\begin{figure}[t]
    \centering
    \includegraphics[scale=0.64]{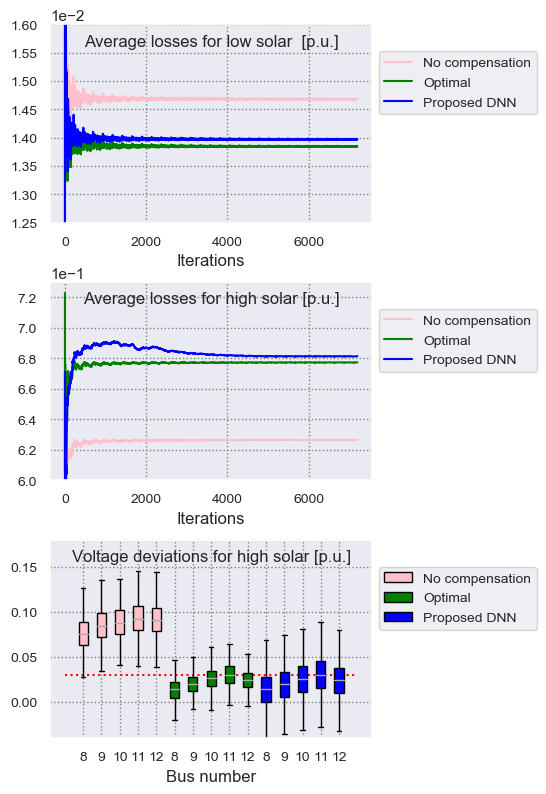}\\
    \hspace*{-7em}\includegraphics[scale=0.75]{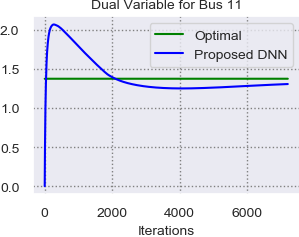}\
    \vspace*{-1em}
    \caption{Training: Average losses under no solar for 12--1~am \emph{(top)} and high solar for 1--2~pm (\emph{(second)}). Voltage excursions for 1--2~pm \emph{(third)}. Dual variable for active constraint on bus 11 for 1--2~pm \emph{(bottom).}}
    \vspace*{-1em}
    \label{fig:results_train}
\end{figure}

We assumed one-hour long control periods. Training scenarios were obtained from the $60$ one-minute data observed over the preceding control period. The original grid scenarios were augmented by adding zero-mean additive white Gaussian noise to generate a total of $K=240$ scenarios. All scenarions were then randomly shuffled. The variance of the additive noise was decided on the basis of training samples observed and was set to $10^{-6}$~pu for low-solar and $10^{-2}$~pu for high-solar hours. DNN $\bpi(\bw;\btheta)$ was trained using Alg.~\ref{alg:1} for $30$ epochs.

The average losses obtained during training are shown in Fig.~\ref{fig:results_train}. The losses under our solution were found to be only slightly superior to those attained by the optimal policy. When inverters operate at unit power factor, significant voltage excursions are observed, while power losses are minimal. Moreover, as demonstrated by the third panel, the proposed scheme attained voltage deviations close to those achieved by the optimal policy. This near-optimal behavior is also shown in the bottom panel presenting the convergence of dual variables for the active constraint on bus 11 during 1:00--2:00~pm.

\begin{figure}[t]
    \centering
    \includegraphics[scale=0.64]{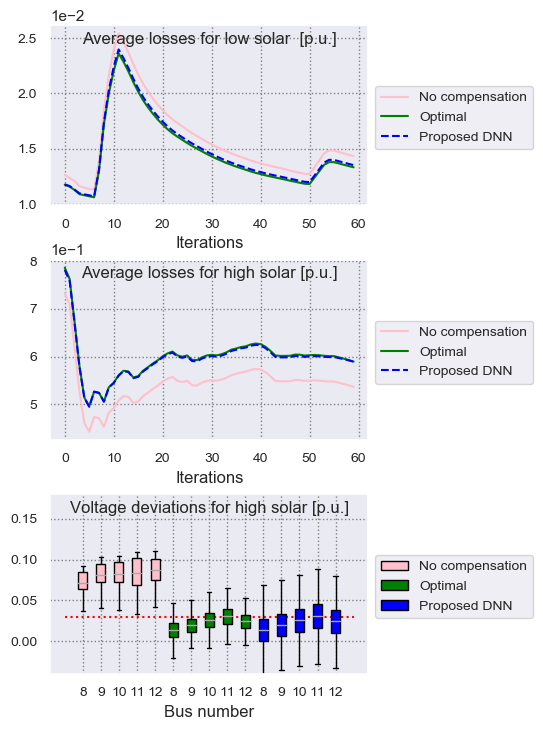}
    \vspace*{-1em}
    \caption{Testing. Average losses for 1--2~am \emph{(top)} and 2--3~pm (\emph{(middle)}); voltage excursions for 2--3~pm \emph{(bottom)}.}
    \vspace*{-1em}
    \label{fig:results_test}
\end{figure}

The DNNs trained over 12:00--1:00~am and 1:00--2:00~pm were tested on the subsequent hours 1:00--2:00~am and 2:00--3:00~pm, respectively. The results are presented in Fig.~\ref{fig:results_test}. The proposed scheme again closely matches the performance of the optimal policy in terms of both minimizing loses and imposing voltage constraints. This is remarkable especially because the optimal policy has access to perfect forecasts and incurs a large real-time communication overhead, while the DNN-based scheme is trained only on historical data and requires only 1 data point to be transmitted in real time.

\section{Conclusions and Ongoing Work}\label{sec:conclusions}
This work has introduced nonlinear control policies for inverter reactive setpoints through a novel two-tier communication-cognizant DNN architecture. The DNN consists of a utility sub-NN and inverter sub-NNs, all jointly trained at the utility at the beginning of every control period, while explicitly incorporating average feeder constraints via primal-dual learning. Upon training, the weights of inverter sub-NNs are downloaded to inverters for real-time implementation. Inverter sub-NNs are driven by local inputs and a control signal broadcast by the utility. Depending on communication specifications, the proposed DNN architecture can accommodate from purely local to centralized and hybrid protocols. Tests on real-world data validate that the suggested methodology is capable of reducing ohmic losses and enforcing feeder constraints with little communication overhead. Furthermore, the \textit{proposed} DNN-based policies were seen to perform comparably to stochastic approximation-based \textit{optimal} policies during both the training and testing phases. 


These promising results set the foundations for relevant generalizations. We are currently working on the following directions: \emph{d1)} Model-free primal-dual learning of DNNs that does not require explicit knowledge of the feeder topology, parameters, and/or precise loading conditions during training; \emph{d2)} Chance-constraint formulations; \emph{d3)} Quantify the performance of the proposed DNN-based approach when compared to the optimal policy; \emph{d4)} incorporating exact AC feeder models; and \emph{d5)} testing on larger feeders to demonstrate scalability.

\balance
\bibliographystyle{IEEEtran}
\bibliography{myabrv,inverters}
\end{document}